\documentclass{notices}
\usepackage[utf8]{inputenc}
\usepackage{amsmath, amssymb}
\usepackage{xcolor}
\usepackage{graphicx}

\newtheorem{theorem}{Theorem}[section]
\newtheorem{example}[theorem]{Example}
\newtheorem{lemma}[theorem]{Lemma}
\newtheorem{condition}[theorem]{Technical condition}
\newtheorem{remark}[theorem]{Remark}

\newcommand{\Br}{\mathrm{Br}}

\title{Cluster algebras in Lie and Knot theory}
\author{Mikhail Gorsky
\affil{University of Vienna. mikhail.gorskii@univie.ac.at. M.~G. received funding from the European Research Council (ERC) under the European Union’s Horizon 2020 research and innovation programme (grant agreement No.~101001159).}
\and
Jos\'e Simental
\affil{Universidad Nacional Aut\'onoma de M\'exico. simental@im.unam.mx. J.S. was partially supported by CONAHCyT Project CF-2023-G-106.}
}
\newcommand{\A}{\mathcal{A}}
\newcommand{\U}{\mathcal{U}}
\newcommand{\Q}{\mathbb{Q}}
\newcommand{\C}{\mathbb{C}}
\newcommand{\T}{\mathbb{T}}
\newcommand{\V}{\mathcal{V}}
\newcommand{\PP}{\mathbb{P}}
\newcommand{\F}{\mathcal{F}}

\newcommand{\R}{\mathbb{R}}
\DeclareMathOperator{\Gr}{Gr}
\DeclareMathOperator{\GL}{GL}
\DeclareMathOperator{\SL}{SL}
\DeclareMathOperator{\Spec}{Spec}
\DeclareMathOperator{\std}{std}
\DeclareMathOperator{\ant}{ant}
\DeclareMathOperator{\HHH}{HHH}
\DeclareMathOperator{\HH}{HH}
\DeclareMathOperator{\gr}{gr}
\date{\vspace{-5ex}}
\begin{document}

\maketitle

%\section{Plan}

%\begin{itemize}
%\item Cluster algebras: few words on history and motivation

%\item Cluster algebras: simplified definition

%\item Example $A_1$ with 1 frozen: algebra

%\item A-Cluster varieties

%\item Example $A_1$ with 1 frozen: geometry
%
%\item Example $A_2$: geometry

%\item Lie theory: short list of varieties with cluster structures

%\item Torus actions and maps to tori 

%\item Cohomology and equivariant cohomology

%\item Braid varieties via Lie theory and configurations of flags

%
%\item Symplectic geometry: links/augmentations
%
%\item Braid closures, Legendrian realizations

%\item Braid barieties as augmentation varieties

%\item "Filling" intuition used to define weaves and construct many clusters

%\item One of the constructions of infinitely many fillings via clusters

%\item Fillings vs. clusters correspondence, works by Korean authors and James

%\item (*) Periphery stuff - depending on which paper gets finished first 

%\item Stratifications

%\item Point count over finite fields, cohomology, etc.

%\item Connections to HOMPFLY and KR homology

\section{Introduction}

Cluster algebras were defined by Sergey Fomin and Andrei Zelevinsky \cite{FZ1} around the beginning of the millenium %\MG{is this grammatically correct? shouldn't it be either ``around the beginning'' or ``towards the beginning''? Jose: No, it was not grammatically corect :)},
with the goal of providing a combinatorial framework for problems related to total positivity and canonical bases in Lie theory. Since then, the theory of cluster algebras has grown to relate to many areas of mathematics, including %\MG{
representation theory and categorification, integrable systems, mathematical physics, higher Teichm\"uller theory, symplectic and algebraic geometry, to name just a few. In this note, we will build from the basics in order to explain a connection between the theory of cluster algebras and that of link invariants, which passes through Lie theory. 

\subsection{Cluster algebras} 
To start, let us define cluster algebras. We will not do it in the greatest possible generality, we refer the reader to \cites{FZ1, FWZ1} for this. The initial data is that of an oriented graph, or a \emph{quiver}, $Q$ with vertex set $Q_0$ that is assumed to be without loops or %\MG{
oriented $2$-cycles. For each vertex $r$ of $Q$, we consider a variable $x_{r}$ and work on the field $\C(x_{r} \mid r \in Q_{0})$. The set $\mathbf{x} := \{x_{r} \mid  r \in Q_{0}\}$ is known as the \emph{initial cluster} and the pair $\Sigma = (Q, \mathbf{x})$ is known as the \emph{initial seed} of the cluster algebra $\A(Q)$. The elements of $\mathbf{x}$ are known as \emph{cluster variables}. 

The main combinatorial input towards defining the cluster algebra is that of a \emph{mutation}. This is associated to a vertex $k$ of the quiver $Q$, and the mutation of a seed $\mu_{k}(\Sigma) = (\mu_{k}(Q), \mu_{k}(\mathbf{x}))$ changes both the quiver $Q$ and the cluster $\mathbf{x}$.  The mutation $\mu_{k}(Q)$ is a new quiver 
%is changes according to  
obtained from $Q$ by the following three-step procedure:
\begin{enumerate}
    \item For each pair of arrows $j \to k \to i$, insert a new arrow $j \to i$. 
    \item Reverse all arrows incident with $k$.
    \item Steps (1) and (2) may have created $2$-cycles. Remove a maximal collection of these. 
\end{enumerate}
%Steps (1)--(3) create a new quiver $\mu_{k}(Q)$. 
The mutation of the cluster $\mu_{k}(\mathbf{x})$ is defined to be $\mu_{k}(\mathbf{x}) = (\mathbf{x} \setminus \{x_{k}\}) \cup \{x'_{k}\}$, where $x'_{k}$ is defined as follows:
\[
x'_{k} = \frac{\prod_{i} x_{i}^{\#\{i \to k\}} + \prod_{i}x_{i}^{\#\{k \to i\}}}{x_{k}}.
\]
The pair $\mu_{k}(\Sigma) := (\mu_{k}(Q), \mu_{k}(\mathbf{x}))$ is known as a \emph{seed}, and $\mu_{k}(\mathbf{x})$ is known as a \emph{cluster}. The elements of $\mu_{k}(\mathbf{x})$ are still called cluster variables. 

One can iterate the mutation procedure indefinitely to obtain an infinite collection of \emph{clusters} and, by definition, the \emph{cluster algebra} $\A(Q) \subseteq \C(x_{r} \mid r \in Q_0)$ is the $\C$-algebra generated by the elements of all these clusters. 

There are certain variations of the above construction that we would like to mention now. We may declare some of the vertices of $Q$ to be unmutable, or \emph{frozen}. We are not allowed to mutate at these vertices. As a consequence, the cluster variable $x_{r}$ belongs to every single cluster, and the cluster algebra $\A$ is a $\C[x_{r} \mid r \; \text{is a frozen vertex}]$-algebra. For these reasons, the frozen cluster variables are also known as \emph{coefficients}. A quiver $Q$ with frozen vertices is also known as an \emph{ice quiver}. It is sometimes convenient to let the cluster variable $x_{r}$ associated to a frozen vertex to be invertible, in which case $\A$ is a $\C[x_{r}^{\pm 1} \mid r \; \text{is a frozen vertex}]$-algebra. This is the convention that we will take for the rest of this article. 

\begin{example}\label{ex:a1 one frozen}
Let us consider the following initial seed
\[
(1 \to {\color{blue} 2}, \{x_1, {\color{blue}x_2}\})
\]
where a blue variable means that it is frozen. We can mutate at vertex $1$ and obtain the new seed
\[
\left(1 \leftarrow {\color{blue} 2}, \left\{x'_1 = \frac{x_{2} + 1}{x_1}, {\color{blue} x_2}\right\}\right).
\]
Mutating at vertex $1$ again, we obtain the initial seed back $(1 \to {\color{blue} 2}, \{x_1, {\color{blue} x_2}\})$. Thus, the associated cluster algebra is $\A = \C[x_{1}, \frac{x_2 + 1}{x_1}, x_{2}^{\pm 1}]$. Note that $x_{2} = x_{1}x'_{1} - 1$, so that we have that $\A$ is the localization $\A = \C[x_{1}, x'_{1}][(x_1x'_1 - 1)^{-1}]$.
\end{example}

Example \ref{ex:a1 one frozen} shows the general phenomenon that \emph{mutation is involutive}. If we mutate a seed twice in the same direction, we go back to the original seed.

\begin{example}\label{ex:a2 one frozen}
For a slightly more complicated example, consider the initial seed $(1 \to 2 \to {\color{blue} 3}, \{x_1, x_2, {\color{blue}x_3}\})$. The reader may verify that by considering all possible mutations we obtain five different clusters: $\{x_1, x_2, {\color{blue} x_3}\}$, $\{\frac{x_2 + 1}{x_1}, x_2, {\color{blue} x_3}\}$, $\{\frac{x_2 + 1}{x_1}, \frac{x_2x_3 + x_3 + x_1}{x_1x_2}, {\color{blue} x_3}\}$, $\{\frac{x_3+x_1}{x_2}, \frac{x_2x_3 + x_3 + x_1}{x_1x_2}, {\color{blue} x_3}\}$, $\{\frac{x_3 + x_1}{x_2}, x_1, {\color{blue} x_3}\}$.
\end{example}

Note that in both examples above, there are only finitely many cluster variables, that is, the corresponding cluster algebra is of \emph{finite cluster type}. The classification of such cluster algebras was one of the early important results in cluster theory.

\begin{theorem}
The cluster algebra associated to an ice quiver $Q$ is of finite cluster type if and only if the quiver obtained after deleting the frozen vertices of $Q$ is mutation equivalent to an orientation of a finite Dynkin diagram.
\end{theorem}

From the definition, it is clear that every cluster variable is a rational function on the initial cluster. Much more is true in the previous example: every cluster variable is actually a \emph{Laurent polynomial} in the initial cluster. This turns out to be a general feature of mutation, that witnesses how fundamental the procedure is.

\begin{theorem}[The Laurent phenomenon]
Let $\A$ be a cluster algebra, and let $\mathbf{x}$ be a cluster of $\A$. Then, any cluster variable of $\A$ can be expressed as a Laurent polynomial in the variables of $\mathbf{x}$.  
\end{theorem}

The Laurent phenomenon was proven by Fomin and Zelevinsky in \cite{FZ1}. They conjectured that, moreover, the coefficients appearing in the Laurent expansion of every cluster variable with respect to a given cluster are all \emph{positive} (as can be readily verified in the previous examples). This positivity conjecture was verified in several cases using several different techniques, most notably involving quiver representations. Lee and Schiffler \cite{LeeSchifflerPositivity} finally settled the positivity conjecture for a very wide class of cluster algebras (including all that we will discuss here), and Gross-Hacking-Keel-Kontsevich \cite{GHKK} proved it in general, using techniques of toric and tropical geometry, as well as mirror symmetry.
%, in \cite{GHKK}.% \MG{and/or tropical geometry (?)} and mirror symmetry, in \cite{GHKK}. 

One way to interpret the Laurent phenomenon is that, for every cluster $\mathbf{x} = \{x_1, \dots, x_r\}$ of $\A$, we have
\[
\A \subseteq \C[x_1^{\pm 1}, \dots, x_{r}^{\pm 1}]
\]
(in fact, $\A[(x_1\cdots x_r)^{-1}] = \C[x_1^{\pm 1}, \dots, x_r^{\pm 1}]$). This leads us to consider the \emph{upper cluster algebra} \cite{BFZ}
\[
\U := \bigcap_{\substack{\mathbf{x} \\ \text{is a cluster}}} \C[x_{1}^{\pm 1}, \dots, x_{r}^{\pm 1}],
\]
so that $\A \subseteq \U$. In nice cases, these two algebras coincide, but this does not always happen. As we will see next, from a geometric point of view it is more natural to consider the algebra $\U$. 

\subsection{Cluster varieties} 
\label{subsection:cluster_varieties}

Yet another, more geometric, point of view on the Laurent phenomenon says that for any cluster $\mathbf{x}$ we have an open torus:
\[
\T_{\mathbf{x}} := \Spec(\C[x_{1}^{\pm 1}, \dots, x_r^{\pm 1}]) \subseteq \Spec(\A)
\]
called a \emph{cluster torus}. By definition, the \emph{cluster manifold} is the union of all cluster tori inside $\Spec(\A)$. Since it is obtained by gluing tori of the same dimension, the cluster manifold is a smooth manifold. Nevertheless, it may not be affine. 

\begin{example}
    Let us go back to Example \ref{ex:a1 one frozen}. Note that, in this case
    \[
    \Spec(\A) = \C^{2} \setminus \{x_1x_2 -1 = 0\}.
    \]
    We have two cluster tori, with coordinates given by $\T_{1} = \{x_1 \neq 0, x_1x_2 - 1\neq 0\}$ and $\T_{2} = \{x_2 \neq 0, x_1x_2 - 1 \neq 0\}$. Thus, the cluster manifold is
    \[
    \T_1 \cup \T_2 = \C^{2} \setminus (\{x_1x_2 - 1 = 0\}\cup \{(0,0)\})
    \]
    which is not an affine variety. 
\end{example}

Note, however, that the algebra of functions which are regular on the entire cluster manifold is precisely the upper cluster algebra $\U$. We define the \emph{cluster variety} 
\[
\V := \Spec(\U),
\]
which is nothing but the affinization of a cluster manifold. We say that an abstract affine algebraic variety $V$ has a cluster structure if its algebra of regular functions $\C[V]$ admits the structure of an upper cluster algebra.

\subsection{Constructing cluster structures}

Given an affine algebraic variety $V$, how to decide whether it admits a cluster structure? Note that, first, one must have candidate for the cluster tori: the variety $V$ must admit a collection of open tori, whose coordinates are given by functions that are regular on the variety $V$. After finding these candidates for cluster tori, the most difficult part is to find a mutation rule that allows us to mutate \emph{every} coordinate of a cluster torus which is not an invertible function on the entire variety $V$ (for these are the frozen variables). This is not an easy ordeal. Nevertheless, we have the following important result, which tells us that it is enough to fix one cluster torus and to be able to mutate in every possible direction at that torus.

\begin{lemma}[The Starfish lemma  \cites{BFZ, FWZ6}]
Let $R$ be a $\C$-algebra that is a Noetherian, normal domain with fraction field $F$. Assume that we are given a seed $\Sigma = (Q, \mathbf{x})$ with $\mathbf{x} \subseteq F$ such that:
\begin{enumerate}
\item The cluster $\mathbf{x}$ consists of elements of $R$.
\item The non-invertible (in $R$) elements of $\mathbf{x}$ are pairwise coprime.
\item For each non-invertible element $x_k \in \mathbf{x}$, mutation at the direction $k$ replaces $x_{k}$ with $x'_{k} \in R$, and $(x_k, x'_k)$ are coprime.
\end{enumerate}
Then, the upper cluster algebra $\U$ associated to $\Sigma$ is contained in $R$. 
\end{lemma}

The geometric proof of this lemma relies heavily on Hartog's lemma, that is valid for normal varieties and that says that a function which is regular outside of a codimension $2$ set must be regular everywhere. The Starfish lemma is very useful when one wants to show that a given algebra has the structure of an (upper) cluster algebra. Note, however, that finding a mutation rule is usually a nontrivial problem, and is usually guided by (sometimes well-hidden) combinatorics.

\subsection{Why?}
Before moving on and giving examples of cluster varieties, we would like to pause and say a few words on what we gain by proving that an algebraic variety $V$ admits a cluster structure. 

In order to explain the first propery below, we have to say that cluster varieties often come in pairs. The cluster varieties we defined in section \ref{subsection:cluster_varieties} are known as \emph{$\mathcal{A}$-cluster varieties}. Fock and Goncharov \cite{FG} introduced a notion of an \emph{$\mathcal{X}$-cluster variety}. This is a Poisson variety associated with the same combinatorial datum of a quiver, but glued from dual algebraic tori. Naturally, the mutation rule for transition maps is different. $\mathcal{A}$- and $\mathcal{X}$-cluster varieties are dual in a sense which has a flavor of mirror symmetry.

In nice cases 
%(for example, the so-called \emph{locally acylic} case, that implies $\A = \U$) \MG{I think this is imprecise: the first property uses the existence of a green-to-red sequence, the third uses some version of the Louise property. If I remember correctly, these are not precisely equivalent. Can we just comment out the part in the brackets?}
one obtains:
\begin{itemize}
\item An explicit basis of the algebra of functions $\C[V]$ with positivity properties, %and that is 
known as the $\vartheta$-basis and parameterized by tropical points of the \emph{dual cluster $\mathcal{X}$-variety} \cites{FG, GHKK}.
% \footnote{In this note, we discuss so-called $\mathcal{A}$-cluster varieties. There is a notion of an $\mathcal{X}$-cluster variety, with a different type of mutation rule, that is dual to an $\mathcal{A}$-cluster variety.}
\item A notion of positivity on the variety $V$.
\item Information about the singular cohomology of $V$. For example, in the locally acylic case, the mixed Hodge structure on $V$ is of mixed Tate type and splits over $\Q$ (in particular, it is a direct sum of pure Hodge structures) \cite{LS}. 
\end{itemize}

Moreover, usually a cluster structure comes with a wealth of combinatorics that allow for explicit computations on the variety $V$.

\section{Cluster varieties in Lie theory}\label{sect:Lie}

Many varieties appearing naturally in Lie theory have been shown to have a cluster structure. Here, we mention only a few. More detailed discussions of these and further examples can be found e.~g. in \cites{FWZ6, GLS}. For concreteness, we restrict ourselves to type $A$, meaning that we will only deal with the Lie groups $\GL(n)$ and $\SL(n)$ of invertible $n \times n$-matrices and $n\times n$ matrices with determinant $1$, respectively.

\begin{enumerate}
\item Let $U \subseteq \SL(n)$ be the subspace of upper triangular matrices with only $1$'s on the diagonal. The \emph{basic affine space} $\SL(n)/U$ is an affine variety whose coordinate ring $\C[\SL(n)]^{U}$ (the ring of invariants) admits a cluster algebra structure, with non-invertible coefficients. The ring $\C[\SL(n)]^{U}$ is generated by the $2^n - 2$ \emph{flag minors}, that is, the determinant of a submatrix occupying the first $1, \dots, k$-rows and $k$ distinct columns ($1 \leq k < n$). Some (but not all) of the clusters are entirely composed of flag minors, and the cluster structure is a consequence of deep determinantal identities. 

\item Let us consider the Grassmannian $\Gr(k, n)$ of $k$-dimensional subspaces in $\C^n$. The \emph{Pl\"ucker embedding} realizes $\Gr(k,n)$ as a closed subvariety of the projective space $\PP^{\binom{n}{k} - 1}$. The affine cone $\widetilde{\Gr}(k, n) \subseteq \C^{\binom{n}{k}}$ over the Grassmannian admits a cluster structure, with non-invertible coefficients. Similarly to the previous case, some (but not all) of the clusters consist entirely of Pl\"ucker coordinates.
\end{enumerate}

We remark that the cluster algebra associated to the Grassmannian $\Gr(k,n)$ is of finite cluster type if and only if $k = 2$, $k = n-2$ or $(k,n) \in \{(3,6), (3,7), (4,7), (3,8), (5,8)\}$.  

\begin{enumerate}
\item[3.] If we take the previous example and consider the corresponding cluster algebra with invertible coefficients, we obtain the coordinate ring of the \emph{maximal positroid cell} $\Pi_{k,n}$ inside the Grassmannian $\Gr(k,n)$. This can be defined as follows. Identify a $k$-dimensional subspace $V \subseteq \C^n$ with a maximal rank $k \times n$ matrix $M_{V}$ whose rows span $V$. Note that this is well-defined up to elementary row operations. Then $\Pi_{k,n}$ is the set consisting of all those subspaces $V$ for which the minors $(M_{V})_{1, 2, \dots, k}, (M_{V})_{2, \dots, k+1}, \dots$ $(M_V)_{n-k+1, \dots, n}$, $M(V)_{n-k+2, \dots, n, 1}, \dots, M(V)_{n, 1, \dots, k-1}$ are all nonzero. In fact, these minors are the frozen variables in (1). 
\item[4.] More general positroid varieties $\Pi \subseteq \Gr(k,n)$, defined by both vanishing and non-vanishing conditions on certain minors of the matrix $M_{V}$. The positroid varieties are affine, smooth and locally closed subsets of $\Pi$. As before, some but not all clusters consist completely of Pl\"ucker coordinates. These clusters can be defined and studied via Postnikov's \emph{plabic  (= planar bicolored) graphs}. 
\item[5.] The space $U \subseteq \SL(n)$ itself has been shown to admit a cluster structure with non-invertible coefficients. Note that, as an algebra, $\C[U]$ is a polynomial ring in $\binom{n}{2}$ variables. 
\end{enumerate}

Cluster algebras are ubiquitous in Lie theory, and have been used to study dual semicanonical bases and total positivity.

\section{Cluster varieties from knot theory}

In the past few years, the theory of cluster algebras has been shown to appear in the study of knots. While there are many instances of this (see e.g. \cites{BMS, LeeSchiffler}),  here we present only one, that originates in the work of Shende-Treumann-Williams-Zaslow \cite{STWZ}, who consider Legendrian links in standard contact $\R^3$. Let us recall that a contact structure on $\R^3$ is a choice of a tangent hyperplane for each point in $\R^3$, satisfying non-integrability conditions. The standard contact structure is choosing the hyperplanes $\ker(ydx - dz)$. A Legendrian link is a smooth, closed $1$-manifold $\Lambda \subseteq \R^3$ whose tangent space at every point is contained in $\ker(ydx - dz)$. In layperson's terms, each component of $\Lambda$ must be given by the image of a smooth embedding $(x(t), y(t), z(t)): S^1 \to \R^3$ satisfying $y(t)x'(t) - z'(t) = 0$ for every $t \in S^1$. Because of this condition, the link $\Lambda \subseteq \R^3$ can be completely recovered from its projection to the $(x,z)$-plane known as the \emph{front projection}. Note that the equation $y(t) = z'(t)/x'(t)$ implies that the projection $\pi(\Lambda) \subseteq \R^2_{xz}$ has no vertical tangencies, and it typically has cusps. 

Given a Legendrian link $\Lambda \subseteq \R^3$, \cite{STWZ} considers the moduli space $\mathcal{M}(\Lambda)$ of microlocal rank $1$ sheaves on $\R^2_{xz}$ microlocally supported on $\pi(\Lambda)$. In \emph{loc. cit.}, the authors construct (not necessarily open) tori in this moduli space $\mathcal{M}(\Lambda)$ and show that, sometimes, one can perform cluster $\mathcal{X}$-mutation on their coordinates, so we have a \emph{partial} $\mathcal{X}$-cluster structure on $\mathcal{M}(\Lambda)$. 
%\MG{This concerns cluster $\mathcal{X}$-transformations and  a partial cluster $\mathcal{X}$-structure. Maybe we should clarify this a bit - but then I'd suggest to slightly expand the content of the footnote in section 1.4 and move it to the main text. The footnote is better  be edited anyway - it's a bit ambiguous what do you mean by ``Here, we discuss'', as it does not mean the sentence to which the footnote is added, but the entire text.} 
Instead of defining $\mathcal{M}(\Lambda)$, we will carefully construct a closely related space, in the case when the Legendrian link $\Lambda = \Lambda(\beta)$ can be constructed from a positive braid $\beta$ in a way we explain next. In order to explain both the link $\Lambda(\beta)$ and its associated algebraic variety, we need to take a detour to braid groups and flag varieties. 

\subsection{Positive braid monoid}

Let us fix a poitive integer $n > 0$. The \emph{positive braid monoid} $\Br^{+}_{n}$ is the monoid generated by elements $\sigma_1, \dots, \sigma_{n-1}$, with relations given by
\begin{equation}\label{eqn:braid relations}
\begin{array}{l}
\sigma_{i}\sigma_{j} = \sigma_{j}\sigma_{i} \qquad  |i - j| > 1,\\
\sigma_{j}\sigma_{j+1}\sigma_{j} = \sigma_{j+1}\sigma_{j}\sigma_{j+1}, \qquad j \leq n-2.
\end{array}
\end{equation}
As its name indicate, the positive braid monoid is nothing but the monoid of positive braids, that is, configurations of $n$ strands where each crossing between the strands is \emph{positive}. Under this point of view, the element $\sigma_{i}$ represents a crossing between the $i$-th and $i+1$-st strands:

\begin{center}
\includegraphics[scale=0.7]{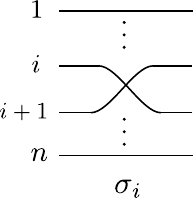}
\end{center}

\noindent Note that, since all the crossings are assumed to be positive, we will not draw over- and under-crossings. The relations \eqref{eqn:braid relations} become purely topological:

\begin{center}
\includegraphics[scale=0.7]{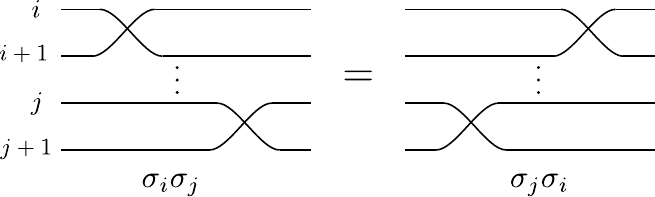}
\end{center}

\begin{center}
\includegraphics[scale=0.6]{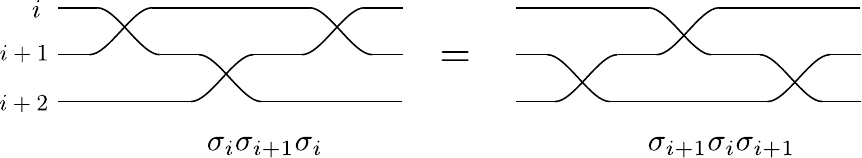}
\end{center}

We say that a braid $\beta \in \Br^{+}_{n}$ is \emph{reduced} if any two strings cross at most once in $\beta$. There is a longest reduced braid $\Delta_{n}$, the one where each pair of strings crosses exactly once in $\beta$. The braid $\Delta_{n}$ is unique, and an expression for it is:
\[
\Delta_n = (\sigma_1\dots \sigma_{n-1})(\sigma_1\dots \sigma_{n-2})\cdots (\sigma_1\sigma_2)\sigma_1.
\]

For example, the element $\Delta_{4} \in \Br_{4}^{+}$ is:
\begin{center}
\includegraphics[scale=0.7]{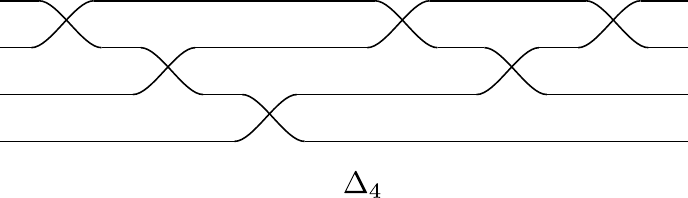}
\end{center}

\begin{center}
\begin{figure}[h!]
\includegraphics[scale=0.5]{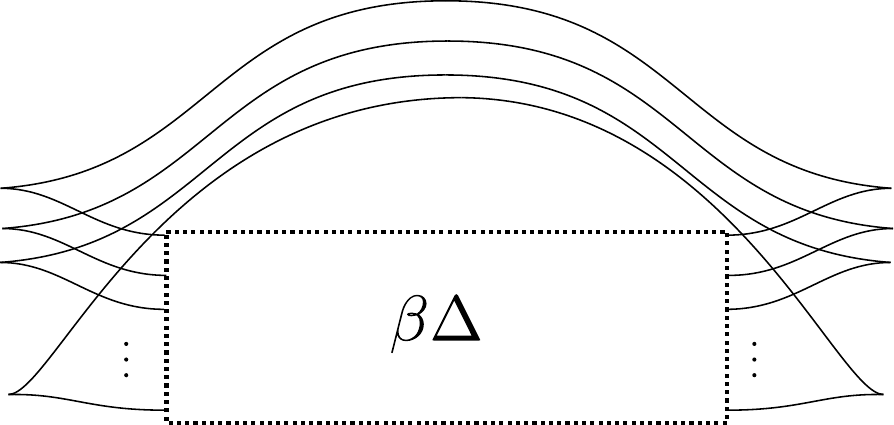}
\caption{Projection to the $xz$-plane of the link $\Lambda(\beta)$. The $y$-coordinate can be recovered by $y = dx/dz$. Note that this implies that all crossings are positive, and that there are no vertical tangencies.}
\label{fig:front}
\end{figure}
\end{center}

Given a braid $\beta \in \Br^{+}_{n}$, we will consider the Legendrian link $\Lambda(\beta)$ whose front projection is as in Figure \ref{fig:front}. It is called the \emph{Legendrian (-1)-closure} of the braid $\beta$.

Finally, we require a technical condition on the braid $\beta$:

\begin{condition}\label{technical}
We will assume that any two strings cross in the braid $\beta$.
\end{condition}%that any two strings cross in it. %\MG{Add some environment and label to cite in section 3.3?}
If we do not require this technical condition, the definition of the variety $X(\beta)$ must be changed accordingly.

\subsection{Flag varieties} 

Having detoured through positive braid monoids, now we need to detour through flag varieties. The flag variety is a projective algebraic variety that appears naturally in Lie representation theory. Consider $n 
> 0$ (this is the same $n$ as in the previous section). We will work on $\C^n$. Let us recall that a \emph{complete flag} $F^{\bullet} = (F^{0} = \{0\} \subseteq F^1 \subseteq \cdots \subseteq F^{n} = \C^n)$ in $\C^n$ is a collection of subspaces $F^{i} \subseteq \C^n$ with $\dim(F^{i}) = i$ and $F_{i} \subseteq F^{i+1}$ for $i = 0, \dots, n-1$. The flag variety $\F_{n}$ is the variety consisting of all complete flags in $\C^n$. It is a smooth projective variety. For example, the flag variety $\F_{2}$ is simply the projective space $\PP^1$. 

The flag variety $\F_n$ can in fact be identified with the quotient space $\SL(n)/B$, where $B \subseteq \SL_{n}$ is the subgroup of upper-triangular matrices. The identification goes as follows. To each matrix $A \in \SL(n)$ we associate the flag $F_{A}^{\bullet}$, where $F_{A}^{i}$ is the subspace spanned by the first $i$ columns of the matrix $A$.  It is easy to see that two matrices $A, B \in \SL(n)$ yield the same flag if and only if $A = UB$, where $U$ is an upper-triangular matrix.

The \emph{standard flag} $F^{\bullet}_{\std}$ is the flag associated to the identity matrix. Equivalently,
\[
F_{\std}^{i} = \langle e_1, \dots, e_{i}\rangle,
\]
where $e_1, \dots e_n$ is the canonical basis of $\C^n$. Similarly, the \emph{antistandard flag} $F^{\bullet}_{\ant}$ is the flag
\[
F_{\ant}^{i} = \langle e_{n}, e_{n-1}, \dots, e_{n-i+1}\rangle.
\]

Two flags $F^{\bullet}$ and $G^{\bullet}$ are said to be \emph{in position} $j \in \{1, \dots, n-1\}$ if
\[
F^{j} \neq G^{j}, \; \text{but} \; F^{i} = G^{i} \; \text{for every} \; i \neq j.
\]
Note that if we fix a flag $F^{\bullet}$, the set of flags that are in position $j$ with respect to $F^{\bullet}$ forms an affine line. Indeed, if $G^{\bullet}$ is in position $j$ with respect to $F^{\bullet}$, all subspaces $G^{i}$ are determined except for $G^{j}$. To choose $G^{j}$ we have to choose a $1$-dimensional subspace in $F^{j+1}/F^{j-1}$ which is different from $F^{j}/F^{j-1}$,  so we have an affine line $\C = \PP^1 \setminus \{\operatorname{pt}\}$ of choices. 

In fact, choosing a matrix $A \in \SL(n)$ such that $F^{\bullet} = F_{A}^{\bullet}$ we have an explicit parametrization of all flags $G^{\bullet}$ which are in position $j$ with respect to $A$. Indeed, these are given by the flags $F^{\bullet}_{AB_{j}(z)}$ for $z \in \C$, where $B_{j}(z)$ is the matrix that looks like the identity everywhere except in the $j$ and $j+1$-st row and column, where it is
\[
B_{j}(z) = \left(\begin{matrix} z & -1 \\ 1 & 0 \end{matrix}\right).
\]
We warn the reader, however, that this parametrization depends on the matrix $A$. If $A'$ is another matrix with $F_{A}^{\bullet} = F_{A'}^{\bullet}$ then it is not the case that $F^{\bullet}_{AB_{j}(z)} = F^{\bullet}_{A'B_{j}(z)}$ for $z \in \C$. Rather, what is true is that for every $z \in \C$ there exists a unique $w \in \C$ with $F^{\bullet}_{AB_{j}(z)} = F^{\bullet}_{A'B_{j}(w)}$. 

\subsection{Braid varieties}
Let us now put things together. Consider a braid $\beta = \sigma_{j_1}\cdots \sigma_{j_{\ell}} \in \Br^{+}_{n}$.  We define the \emph{braid variety}
\[
X(\beta) := \{F_{0}^{\bullet}, F_{1}^{\bullet}, \dots, F_{\ell}^{\bullet}\}
\]
to be the space consisting of sequences of flags satisfying the following conditions:
\begin{enumerate}
\item $F_{0}^{\bullet} = F_{\std}^{\bullet}$.
\item $F_{\ell}^{\bullet} = F_{\ant}^{\bullet}$.
\item For every $i = 1, \dots, \ell$, the flags $F^{\bullet}_{i-1}$ and $F^{\bullet}_{i}$ are in position $j_{i}$. 
\end{enumerate}

Note that the braid variety is \emph{affine}. Indeed, since $F_{0}^{\bullet} = F_{\std}^{\bullet}$, one may describe the braid variety as follows:
\[
X(\beta) = \{(z_1, \dots, z_{\ell}) \in \C^{\ell} \mid F^{\bullet}_{B_{j_1}(z_1)\cdots B_{j_{\ell}}(z_{\ell})} = F^{\bullet}_{\ant}\}.
\]
The condition $F^{\bullet}_{B_{j_1}(z_1)\cdots B_{j_{\ell}}(z_{\ell})} = F^{\bullet}_{\ant}$ can be more explicitly described as
\begin{equation}\label{eqn:def braid variety}
w_{0}^{n}B_{j_1}(z_1)\cdots B_{j_{\ell}}(z_{\ell}) \; \text{is upper-triangular},
\end{equation}
where $w_0^{n}$ is the $n \times n$-matrix given by $(w_{0}^{n})_{i,j} = \delta_{i, n-j+1}$, that is, $w_{0}^{n}$ has $1$'s on the main anti-diagonal and $0$'s everywhere else. Note that \eqref{eqn:def braid variety} is equivalent to the vanishing of $\binom{n}{2}$ entries of the matrix $B_{j_1}(z_1)\cdots B_{j_{\ell}}(z_{\ell})$, and the entries of the latter matrix are polynomials in $z_1, \dots, z_{\ell}$. Thus, we indeed have an affine variety. 

It is instructive to look at the case $n = 2$. In this case, $\beta = \sigma_1^{\ell}$, and the flags $F_{i}^{\bullet}$ are simply elements of $\mathbb{P}^{1} = \C \cup \{\infty\}$. The standard flag is identified with $0 \in \C \subseteq \mathbb{P}^1$, and the antistandard flag with $\infty \in \mathbb{P}^1$. Thus, $X(\sigma^{\ell})$ consists of collections of points
\[
(0 = x_0, x_1, \dots, x_{\ell} = \infty) \in (\mathbb{P}^{1})^{\ell + 1}
\]
satisfying $x_{i} \neq x_{i+1}$ for every $i = 0, \dots, \ell - 1$. Thus, $X(\sigma_1^{2}) = \mathbb{P}^{1}\setminus\{0, \infty\} = \C^{\times}$. 

Note that, a priori, the variety $X(\beta)$ depends not only on the brad $\beta$, but also on the chosen expression $\beta = \sigma_{j_1}\cdots \sigma_{j_{\ell}}$ for it. In fact, up to a canonical isomorphism, the variety $X(\beta)$ does indeed depend \emph{only} on the braid $\beta$. The name \emph{braid variety} is chosen to emphasize this. %Let us explain why the name \emph{braid} variety. If $\bj = (j_1, \dots,j_{i-1}, j_{i}, j_{i+1}, \dots, j_{\ell})$ and $\bj' = (j_{1}, \dots, j_{i-1}, j_{i+1}, j_{i}, \dots, j_{\ell})$, then $X(\bj) \cong X(\bj')$ provided $|j_i - j_{i+1}| > 1$. Moreover, if $\bj = (j_1, \dots, k, k+1, k, \dots, j_{\ell})$ and $\bj' = (j_1, \dots, k+1, k, k+1, \dots, j_{\ell})$ then we also have $X(\bj) = X(\bj')$. On the other hand, the \emph{positive braid monoid} $\Br_{n}^{+}$ is the monoid defined by generators $\sigma_1, \dots, \sigma_{n-1}$ and relations:

Note that $X(\beta)$ may be empty. For example, if $n = 3$ $X(\sigma_1) = \emptyset$: there is simply no space to start with $F^{\bullet}_{\std}$ and end with $F^{\bullet}_{\ant}$. This is where the technical condition \ref{technical} on $\beta$ is needed: the variety $X(\beta)$ is nonempty if and only if $\beta$ satisfies this technical condition. 

Thanks to results of Brou\'e-Deligne-Michel, the variety $\mathcal{M}^{1}(\Lambda(\beta))$ is closely related to $X(\beta)$. Braid varieties have, however, appeared in many other places. Most notably, they appeared in the work of Mellit \cite{Mellit}, who stratifies character varieties using braid varieties in order to prove the curious Lefschetz property for the former. A compactification of braid varieties, in which two consecutive flags are not required to be different, appeared in the work of Escobar \cite{Escobar} under the name of \emph{brick varieties}. These are general fibers of the moment maps on certain symplectic varieites with Hamiltonian torus actions. 
%used the compactification to explicitly construct brick polytopes. 
We remark, however, that such a compactification depends on a chosen word for $\beta$ and not on $\beta$ itself. 

\section{Weaves} \label{section:weaves}

The braid variety $X(\beta)$ admits a cluster structure, where the cluster tori are parametrized by objects called \emph{(Demazure) weaves}. A weave is a certain colored graph, that is similar to the graphical calculus appearing in Soergel theory. In the setting of Legendrian geometry, they were recently introduced by Casals and Zaslow \cite{CasalsZaslow}, who use weaves in order to construct exact Lagrangian surfaces in $\R^4$ %\MG{in $\R^5?$} \JS{I think it's $\R^4$. Legendrian in $\R^5$ whose projection to $\R^4$ is Lagrangian} 
whose boundary is precisely the Legendrian link $\Lambda(\beta) \subseteq \mathbb{R}^3 \subseteq \mathbb{S}^{3}$, that is, an \emph{exact Lagrangian filling} of $\Lambda(\beta)$. We will, however, not take this point of view.

In order to motivate weaves, let us take $\beta = \sigma_{j_1}\cdots \sigma_{j_{\ell}}$, and we picture an element of the braid variety $X(\beta)$ as follows:

\begin{center}
\includegraphics[scale=0.5]{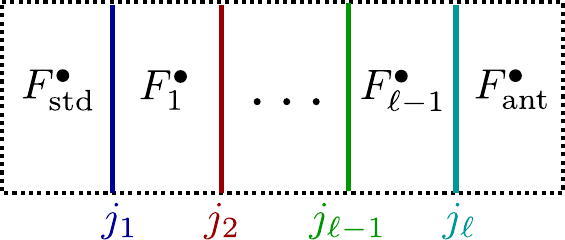}
\end{center}

Here, the vertical lines are colored with colors $1, \dots, n-1$. The flag on the leftmost region is $F_{\std}^{\bullet}$, and the flag on the rightmost is $F_{\ant}^{\bullet}$. If two regions are separated by a line of color $i$, then the corresponding flags are in position $i$, that is, they differ precisely in the $i$-th subspace.

Now suppose we have the following configuration:

\begin{center}
\includegraphics[scale=0.5]{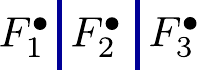}
\end{center}

The flags $F_{1}^{\bullet}$ and $F_{3}^{\bullet}$ are either equal, or they differ in precisely the $i$-th subspace. It is natural to picture these two possibilities using the following diagrams.

\begin{center}
    \includegraphics[scale=0.5]{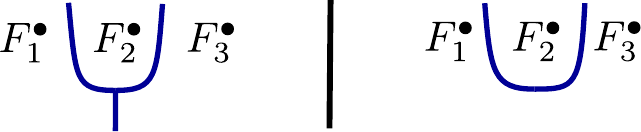}
\end{center}

\begin{center}
\begin{figure}
\centering
\includegraphics[scale=0.5]{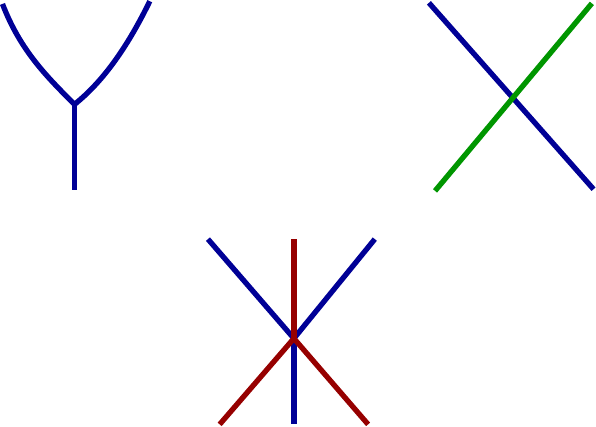}
\caption{The types of vertices in the interior of the rectangle $R$ of the definition of a weave. Note that the edges adjacent to a trivalent vertex all have the same color; the edges adjacent to a tetralent vertex are of $2$ distant colors; and the edges adjacent to an hexavalent vertex are of neighboring colors.}
\label{fig:types of vertices}
\end{figure}
\end{center}

Note that the condition $F_{1}^{\bullet} \neq F_{3}^{\bullet}$ is \emph{open}, and the complement $F_{1}^{\bullet} = F_{3}^{\bullet}$ is closed. A Demazure weave $\mathfrak{w}$ on $\beta$ then is a graph in a rectangle $R$, whose edges are colored $1, \dots, n-1$ and whose vertices are one of four types:
\begin{enumerate}
    \item Univalent vertices, located only on the top and bottom sides of $R$. Moreover, the edges adjacent to the vertices on the top side spell a braid word for $\beta$, and the word spelt by the colors of the edges adjacent to the vertices on the bottom side is reduced.
    \item Trivalent vertices, as pictured in Figure \ref{fig:types of vertices}.
    \item Tetravalent vertices, as pictured in Figure \ref{fig:types of vertices}.
    \item Hexavalent vertices, as pictured in Figure \ref{fig:types of vertices}. 
\end{enumerate}

\begin{figure}
    \centering
    \includegraphics[scale=1]{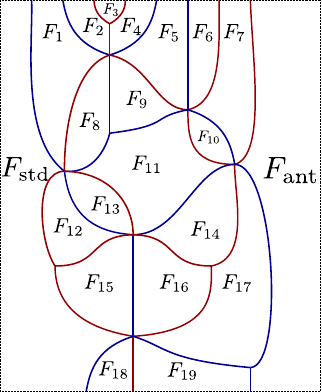}
    \caption{A weave on $\beta = \sigma_1^{2}\sigma_{2}^{2}\sigma_{1}^{2}\sigma_{2}^{2}$ and its flag moduli. All the flags are determined by $F_{\std}, F_{1}, \dots, F_{7}$ and $F_{\ant}$.}
    \label{fig:weave moduli}
\end{figure}

An example of a Demazure weave is given on Figure \ref{fig:weave moduli}, taken from \cite{CGGLSS}.

Given a weave $\mathfrak{w}$ on $\beta$, we consider the moduli space $X(\mathfrak{w})$ of all configurations of flags $F^{\bullet}_{C}$, one per connected component $C$ of $R \setminus \mathfrak{w}$, satisfying the following conditions. 

\begin{itemize}
\item The flag labeling the region bordering the left side of the rectangle is $F^{\bullet}_{\std}$.
\item The flag labeling the region bordering the right side of the rectangle is $F^{\bullet}_{\ant}$.
\item If two regions are separated by an edge of color $i$, then the corresponding flags are in position $i$.
\end{itemize}

By definition, it is clear that looking at the flags labeling the regions bordering the \emph{top} side of the rectangle, we obtain an element of $X(\beta)$. In fact, these flags propagate to fill the entire rectangle, and we see that actually $X(\mathfrak{w}) \subseteq X(\beta)$. The \emph{trivalent vertices} impose conditions on the element of $X(\beta)$ and, moreover, $X(\mathfrak{w})$ is an open torus inside $X(\beta)$. See Figure \ref{fig:weave moduli}.

The tori $X(\mathfrak{w})$ are cluster tori in a cluster structure on $X(\beta)$. Thus, each weave determines an ice quiver $Q_{\mathfrak{w}}$ and a collection of regular functions on $X(\beta)$, the cluster variables. The determination of the cluster variables is technical and we will not explain it here, it suffices to say that  they are given by explicit polynomials in the coordinates $z_1, \dots, z_{\ell}$ of $X(\beta)$. 

\begin{figure}
    \centering
    \includegraphics[scale=0.85]{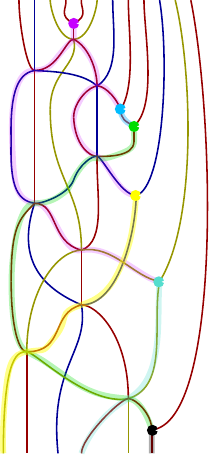}
    \caption{A weave on $\beta = \sigma_2\sigma_1\sigma_3\sigma_2^{2}\sigma_3\sigma_1\sigma_2^{2}\sigma_1\sigma_3\sigma_2$ (we omit the rectangle). The paths starting at every trivalent vertex are indicated on the figure.}
    \label{fig:longweave}
\end{figure}

The vertices of the ice quiver $Q_{\mathfrak{w}}$ are in correspondence with the trivalent vertices of the weave $\mathfrak{w}$. To determine the arrows, we define a collection of $d$ positive linear combinations of
paths on the weave itself, 
%\MG{should we clarify that weights on edges can be larger than one?}, starting at the corresponding trivalent vertex and 
obeying certain rules known as \emph{Lusztig's tropical rules}, see Figure \ref{fig:longweave} for an example where each linear combination consists of a single path. Here $d$ is the number of trivalent vertices in $\mathfrak{w}$, which equals the dimension of $X(\beta)$. The arrows are then given by signed intersections between these linear combinations of paths. It is interesting to note that the paths define homology cycles on the Lagrangian surface associated to the weave $\mathfrak{w}$ by Casals--Zaslow \cite{CasalsZaslow}, and these intersections are topological intersection numbers of cycles. %\MG{``intersection numbers''?} 

For example, the quiver corresponding to the weave in Figure \ref{fig:longweave} is

\begin{center}
    \includegraphics[scale=0.8]{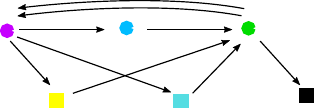}.
\end{center}

\noindent Here the frozen vertices are depicted as squares -- note that they correspond to those trivalent vertices whose corresponding paths reach the bottom of the weave. Note also that there are two intersections between the fuchsia and the green paths in the weave in Figure \ref{fig:longweave}, so we have two arrows between the corresponding vertices in the quiver. 

There is a notion of \emph{weave mutation}, as shown in the following figure:

\begin{center}
    \includegraphics[scale=0.5]{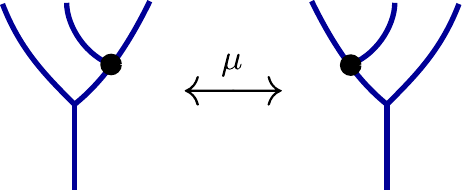}
\end{center}

\noindent This weave mutation corresponds to the quiver mutation at the vertex associated to the indicated trivalent vertex of the weave. The operation $\mu$ is clearly involutive on weaves, although it does not allow us to mutate at every vertex of the weave. More general weave mutations are described in \cite{CasalsZaslow}.%[Casals-Weng] \MG{In [Casals-Zaslow]?}.

The cluster algebra $\A$ associated to the braid variety $X(\beta)$ has very nice properties:
\begin{itemize}
    \item It is \emph{locally acyclic}, that is, $X(\beta)$ can be covered by cluster varieties associated to cluster algebras whose defining quiver is acyclic, \cite{Muller}.
    \item It is equal to its own upper cluster algebra.
    \item It is really full rank, \cites{CGGLSS, GLSBS}. 
    \item The variety $X(\beta)$ also admits the structure of a cluster $\mathcal{X}$-variety \cite{FG}. Moreover, $X(\beta)$ is its own cluster dual.
    \item It admits a $\vartheta$-basis, \cite{GHKK}.
\end{itemize}

%\MG{Almost none of these have been defined. It is mentioned earlier that local acyclicity implies A = U, so we can switch the order of the first two points. Also, add X-structures to the list.}

Let us notice that, although the cluster structure and cluster variables in $X(\beta)$ can be obtained in a purely combinatorial way, they should also be, in principle, obtained using symplectogeometric tools. This has been explicitly realized by Casals and Weng in \cite{CasalsWeng} for weaves coming from \emph{grid plabic graphs}, and it is an interesting problem to do it for general weaves. %\MG{Yet another construction of cluster structures on braid varieties appeared in \cite{GLSBS} and uses a 3-dimensional analogue of plabic graphs for combinatorial parameterisation of some of the seeds.}

We note that yet another construction of cluster structures on braid varieties appeared in \cite{GLSBS} and uses a 3-dimensional analogue of plabic graphs for combinatorial parameterisation of some of the seeds. The connection between the 3D plabic gaphs of \cite{GLSBS} and the weaves we define here remains to be elucidated. 

We remark that many of the varieties mentioned in Section \ref{sect:Lie} admit a realization as braid varieties. In particular, each open positroid stratum in $\Gr(k,n)$ can be realized as a braid variety for a braid with $n$-strands and also (up to a torus factor) as a braid variety for a braid with $k$ strands. 

\section{Applications}

\subsection{Cluster structures on Richardson varieties}
Let $w \in S_{n}$. We say that two flags $F^{\bullet}$ and $G^{\bullet}$ are in position $w$, and write $F^{\bullet} \buildrel w \over \rightarrow G^{\bullet}$ if, for every $i, j \in \{1, \dots, n\}$:
\[
\dim(F^{i}\cap G^{j}) = \# w([1, \dots, i])\cap[1, \dots, j],
\]
for example, $F^{\bullet} \buildrel s_{i} \over \longrightarrow G^{\bullet}$ if and only if $F^{\bullet}$ and $G^{\bullet}$ are in position $i$. The \emph{Schubert cell} associated to $w$ is
\[
C_{w} = \{F^{\bullet} \in \F_{n} \mid F^{\bullet}_{\std} \buildrel w \over \rightarrow F^{\bullet}\}.
\]
This is known to be an affine space of dimension $\ell(w) := \{(i, j) \in \{1, \dots, n\} \mid i < j \; \text{and} \; w(j) < w(i)\}$. Similarly, the \emph{opposite Schubert cell}
\[
C^{w} := \{F^{\bullet} \in \F_{n} \mid F^{\bullet} \buildrel w_{0}w \over \rightarrow F^{\bullet}_{\ant}\}
\]
is an affine space of dimension $\binom{n}{2} - \ell(w)$. Note the presence of $w_0$ in the definition of $C^{w}$: this is the longest permutation in $S_{n}$, defined by $w_0(i) = n+1-i$. For $v, w \in S_{n}$ the \emph{open Richardson variety} $R(w,v)$ is defined to be
\[
R(w,v) := C_{w} \cap C^{v}.
\]
Open Richardson varieties appear naturally in Lie theory, since the cohomology of Richardson varieties naturally computes extensions between Verma modules. These are certain infinite-dimensional representations of $\mathfrak{sl}(n)$ forming the standard modules in the principal block $\mathcal{O}_{0}$ of the Berstein-Gelfand-Gelfand category $\mathcal{O}$, a highest-weight category naturally associated to $\mathfrak{sl}(n)$.

Richardson varieties can be naturally realized as braid varieties, and thus we obtain a cluster structure on Richardson varieties. To realize $R(w, v)$ as a braid variety, consider reduced (i.e. minimal length) expressions
\[
w = s_{i_1}\cdots s_{i_{k}}, \qquad v^{-1}w_{0} = s_{j_1}\cdots s_{j_{\ell}},
\]
then
\[
R(w, v) \cong X(\sigma_{i_{1}}\cdots \sigma_{i_{k}}\sigma_{j_{1}}\cdots \sigma_{j_{\ell}}).
\]
It is known that the Richardson variety is empty unless $v \leq w$ in Bruhat order -- which translates to the Technical Condition \ref{technical} on the braid variety side.

Richardson varieties were conjectured to have cluster structures by Leclerc in \cite{Leclerc}. It is known that Schubert cells have cluster structures. We remark that positroid varieties are a special case of Richardson varieties, obtained when the permutation $w$ is $k$-Grassmannian, that is:
\[
w^{-1}(1) < \cdots < w^{-1}(k), \quad w^{-1}(k+1) < \cdots < w^{-1}(n),
\]
in fact, when $w$ is $k$-Grassmannian, the restriction of the natural projection $\pi_{k}: \F_{n} \to \Gr(k,n)$ to the Schubert cell $C_{w}$ (and thus to any Richardson variety of the form $R(w,v)$) is an isomorphism. The relationship between the cluster structures on positroid varieties obtained via weaves and via plabic graphs
was described very recently in \cite{CLSBW}. 
%remains to be described. 

\subsection{Infinitely many Lagrangian fillings}

Let $\Lambda \subseteq \R^3 \subseteq \mathbb{S}^{3}$ be a Legendrian link. An \emph{exact Lagrangian filling} of $\Lambda$ is an exact Lagrangian surface $S \subseteq \mathbb{B}^{4} \subseteq \mathbb{R}^{4}$ whose boundary $\partial S$ coincides with $\Lambda$. Any two such fillings are known to be \emph{smoothly} isotopic. For the purposes of symplectic geometry, the correct notion of isotopy is that of \emph{Hamiltonian} isotopy, i.e., an isotopy given by a family of Hamiltonian vector fields on the symplectic $4$-ball $\mathbb{B}^{4}$. 

The classification of Lagrangian fillings up to Hamiltonian isotopy is one of the central problems in contact and symplectic geometry. A complete classification of Lagrangian fillings is known only for the standard Legendrian unknot. Interpreted in the way of Casals-Zaslow, a weave on $\beta$ gives an exact Lagrangian filling of the link $\Lambda(\beta)$. Moreover, if the Lagrangian fillings given by the weaves $\mathfrak{w}, \mathfrak{w}'$ are Hamiltonian isotopic, then the cluster tori $X(\mathfrak{w})$ and $X(\mathfrak{w}')$ coincide.

While the correspondence between fillings and cluster tori is still not entirely well understood, Casals and Gao \cite{CG} used this circle of ideas in order to give the first examples of \emph{infinitely many} (not Hamiltonian isotopic) Lagrangian fillings of Legendrian links. More precisely, they show that a Legendrian torus link
\[
\Lambda(n,m) := \Lambda((\sigma_{1}\cdots \sigma_{n-1})^{m}\Delta_{n}), \qquad n \leq m
\]
admits infinitely many Lagrangian fillings unless $(n,m) \in \{(2,m), (3,3), (3,4), (3,5)\}$. These exceptional cases correspond to Grassmannians which are of finite cluster type. 

It should be emphasized that, while \emph{loc. cit.} does not explicitly use cluster algebras, the techniques are heavily inspired by cluster techniques, in particular, they use actions of cluster modular groups in order to produce infinitely many fillings. These techniques were later extended by Gao, Shen and Weng in \cite{GSW} who in fact prove that, for braids of the form $\beta = \gamma\Delta$, the Legendrian link $\Lambda(\beta)$ admits infinitely many Lagrangian fillings unless the cluster algebra $\C[X(\beta)]$ is of finite cluster type. Moreover, they show that for such braid words $\beta = \gamma\Delta$, $\C[X(\beta)]$ is of finite cluster type if and only if $\Lambda(\gamma)$ is Legendrian isotopic to a split union of unknots and the following types of Legendrian links:
\begin{enumerate}
\item[($A_n$)] $\Lambda(\sigma_{1}^{n+1}\Delta_2)$.
\item[($D_n$)] $\Lambda(\sigma_1^{n-2}\sigma_2\sigma_1^{2}\sigma_2\Delta_3)$.
\item[($E_6$)] $\Lambda(\sigma_1^{3}\sigma_2\sigma_1^{3}\sigma_2\Delta_3)$.
\item[($E_7$)] $\Lambda(\sigma_1^{4}\sigma_2\sigma_1^{3}\sigma_2\Delta_3)$.
\item[($E_8$)] $\Lambda(\sigma_1^{5}\sigma_2\sigma_1^{3}\sigma_2\Delta_3)$.
\end{enumerate}
The labels, of course, correspond to the Dynkin type of the mutable part of the corresponding quiver. Finally, we remark that, even if the cluster algebra $\C[X(\beta)]$ is of finite cluster type, it is \emph{not known} whether the link $\Lambda(\beta)$ admits infinitely many Lagrangian fillings. 

\begin{remark}
    Very recently, Casals and Gao \cite{CG2023} have shown that, for braids of the form $\beta = \gamma\Delta_{n}$, with $\gamma \in \Br^{+}_{n}$, \emph{every} cluster comes from a filling in a precise way, see \cite{CG2023}. It is not known whether the same cluster can correspond to more than one filling.
\end{remark}

\subsection{Khovanov-Rozansky homology}

In this last section, we forget about Legendrian structures and work only with \emph{smooth} links. Khovanov-Rozansky homology is a powerful invariant of smooth links that is, nevertheless, notoriously hard to compute. To a link $L$, Khovanov-Rozansky homology associates a triply graded vector space $\HHH(L)$ that is an invariant of $L$, and  whose Poincar\'e polynomial recovers, under various specializations, more classical link invariants such as the Jones and HOMFLY-PT polynomials. 

The construction of $\HHH(L)$ is highly technical. First, one must present $L$ as the $0$-framed closure of a braid $\beta$, $L = L(\beta)$, that is, as follows:
\begin{center}
\includegraphics[scale=0.5]{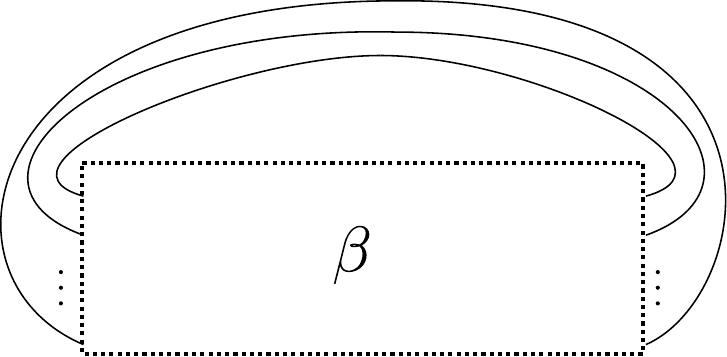}
\end{center}
Assuming $\beta$ has $n$ strands, one then constructs a complex of graded $\C[x_1, \dots, x_n]$-bimodules, known as the \emph{Rouquier complex} $T_{\beta}$ of $\beta$. Then, one takes Hochschild homology of the complex $T_{\beta}$, to obtain a new complex of graded $\C[x_1, \dots, x_n]$-bimodules $\HH(T_{\beta})$. Then, $\HHH(L)$ is defined to be the homology of $\HH(T_{\beta})$. The three gradings on $\HHH(L)$ come, essentially, from the internal grading of $\C[x_1, \dots, x_n]$-bimodules, the homological grading on the complex $T_{\beta}$, and the Hochschild grading on $\HH(T_{\beta})$. These degrees are denoted by $q, t$ and $a$, respectively. 

Due to the highly technical nature of $\HHH(L)$, geometric models have been proposed for it. One of them, going back to work of Webster-Williamson \cite{WW} and recently established by Trinh \cite{Trinh}, is closely related to braid varieties. First, for a positive braid $\beta$ on $n$ strands, we consider the variety $\widetilde{X}(\beta)$, whose definition closely models that of the braid variety $X(\beta)$ except that we require both the first and last flags to be the standard flag: $F_{0}^{\bullet} = F_{\ell}^{\bullet} = F_{\std}^{\bullet}$. Let us assume that $\beta = \gamma\Delta$ for a positive braid $\gamma$. Then it is possible to show that
\[
\widetilde{X}(\beta) = X(\gamma) \times \C^{\binom{n}{2}},
\]
so that the varieties $\widetilde{X}(\beta)$ and $X(\gamma)$ are closely related. The torus $T = (\C^{\times})^{n}$ acts on the smooth variety $\widetilde{X}(\beta)$. The $T$-equivariant Borel-Moore homology of $\widetilde{X}(\beta)$ has a nontrivial \emph{weight filtration} $\mathrm{W}$ and Trinh \cite{Trinh} showed that
\[
\begin{array}{rl}
\gr_{\mathrm{W}}H^{T}_{*, BM}(\widetilde{X}(\beta))  \cong & 
\HHH^{a = n}(L(\beta))\\ \cong &  \HHH^{a = 0}(L(\gamma\Delta^{-1}))
\end{array}
\]
so that the homology of the braid variety $X(\gamma)$ recovers part of the Khovanov-Rozansky homology of the link $L(\gamma\Delta^{-1})$. In many cases, the action of the torus $T$ coincides with the action of the so-called \emph{torus of cluster automorphisms} on the cluster variety $X(\gamma)$. Thus, cluster-theoretic arguments can be used to compute (part of) the Kovanov-Rozansky homology of certain links. This has been done by Galashin and Lam \cite{GL} for so-called \emph{positroid links}, and remains to be done in more generality.

\end{document}